The mod-$p$ cohomology rings of some $p$-groups.
I.J. Leary*
School of Mathematical Sciences,
Queen Mary and Westfield College,
Mile End Road,
London.
E1 4NS

**Introduction.**

Throughout this paper $p$ denotes an odd prime. The groups we shall consider are central extensions of a cyclic subgroup by $C_p \oplus C_p$, and may be presented as

$$P(n) = \langle A, B, C \mid A^p = B^p = C^{p^{n-2}} = [A,C] = [B,C] = 1, \quad [A,B] = C^{p^{n-3}} \rangle.$$

The group $P(n)$ is defined for each $n \geqslant 3$, and has order $p^n$. The group $P(3)$ is the nonabelian group of order $p^3$ and exponent $p$. The mod-$p$ cohomology rings of the other groups of order $p^3$ have been known for some time (see [14] or [5] for that of the nonabelian metacyclic group), and so this paper completes the calculation of the mod-$p$ cohomology rings of the groups of order $p^3$. The corresponding calculations for integral cohomology were completed by Lewis [11]. In the case when $p = 3$, Milgram and Tezuka have determined the mod-$p$ cohomology of $P(3)$ by showing that it is detected by the restrictions to proper subgroups [13]. For $p \geqslant 5$ the cohomology of $P(3)$ is not detected by restrictions to proper subgroups, so the method of [13] cannot be used. There is some overlap between this paper and the paper of Benson and Carlson on the cohomology of extraspecial groups [2].

We calculate $H^*(\mathrm{B}P(n); \mathbb{F}_p)$ using a similar method to our calculation of the integral cohomology of the same groups [10]. First, we let $\widetilde{P}$ be the unique nonabelian Lie group whose underlying topological space consists of $p^2$ circles. $\widetilde{P}$ may be 'presented' as follows, where we consider $\mathrm{S}^1$ as a subgroup of $\mathbb{C}$.

$$\widetilde{P} = \langle \mathrm{S}^1, A, B \mid A^p = B^p = 1, \ \mathrm{S}^1 \text{ central}, [A,B] = \exp(2\pi i/p) \rangle.$$

The subgroup of $\widetilde{P}$ generated by $A$, $B$, and $C = \exp(2\pi i/p^{n-2})$ is a normal subgroup isomorphic to $P(n)$, with quotient isomorphic to $\mathrm{S}^1$, and from now on we shall regard $P(n)$ as a subgroup of $\widetilde{P}$, using this fixed embedding. It follows that $\mathrm{B}P(n)$ is a principal $\mathrm{S}^1$-bundle over $\mathrm{B}\widetilde{P}$. We determine $H^*(\mathrm{B}\widetilde{P}; \mathbb{F}_p)$, and then use the two row spectral sequence (or Gysin sequence) for the fibration

$$\mathrm{S}^1 \longrightarrow \mathrm{B}P(n) \longrightarrow \mathrm{B}\widetilde{P}$$

to find $H^*(\mathrm{B}P(n); \mathbb{F}_p)$. This method, which may be applied to any group with nontrivial centre, was suggested independently by P.H. Kropholler and J. Huebschmann [6], [7].

* Supported by SERC post-doctoral fellowship B90 RFH 8960.

**Calculations.**

We begin by calculating $H^*(\mathrm{B}\widetilde{P};\mathbb{F}_p)$. This may be achieved by using the spectral sequence for $\widetilde{P}$ as an extension of $\mathrm{S}^1$ by $C_p \oplus C_p$, which has only two non-zero differentials, $d_3$ and $d_{2p-1}$. This approach was taken in the author's thesis [9], and in work of B. Moselle [15], who found $H^*(\mathrm{B}\widetilde{P};\mathbb{F}_p)$ while investigating extraspecial $p$-groups. It is also possible to calculate $H^*(\mathrm{B}\widetilde{P};\mathbb{F}_p)$ directly from $H^*(\mathrm{B}\widetilde{P};\mathbb{Z})$, and this is the method we shall employ. The following statement is taken from theorem 2 of [10].

**Theorem 1.** *Let $\widetilde{P}$ be the group defined above. Then $H^*(\widetilde{P};\mathbb{Z})$ is generated by elements $\alpha, \beta, \chi_1, \ldots, \chi_{p-1}, \zeta$, with*

$$\deg(\alpha) = \deg(\beta) = 2 \quad \deg(\chi_i) = 2i \quad \deg(\zeta) = 2p,$$

*subject to the following relations:*

$$p\alpha = p\beta = 0$$

$$\alpha^p \beta = \beta^p \alpha$$

$$\alpha\chi_i = \begin{cases} 0 \\ -\alpha^p \end{cases} \quad \beta\chi_i = \begin{cases} 0 & \text{for } i < p-1 \\ -\beta^p & \text{for } i = p-1 \end{cases}$$

$$\chi_i \chi_j = \begin{cases} p\chi_{i+j} & i+j < p \\ p^2\zeta & i+j = p \\ p\zeta\chi_{i+j-p} & p < i+j < 2p-2 \\ p\zeta\chi_{p-2} + \alpha^{2p-2} + \beta^{2p-2} - \alpha^{p-1}\beta^{p-1} & i = j = p-1 \end{cases}$$

*An endomorphism of $\widetilde{P}$ that restricts to $\mathrm{S}^1$ as $z \mapsto z^j$ sends $\chi_i$ to $j^i \chi_i$ and $\zeta$ to $j^p \zeta$. The effect of an automorphism on $\alpha$, $\beta$ may be determined from their definition. Considered as elements of $\mathrm{Hom}(\widetilde{P}, \mathrm{S}^1)$, $\alpha$ and $\beta$ may be defined by $\alpha(A^a B^b z) = \eta^a$, $\beta(A^a B^b z) = \eta^b$, where $\eta = \exp 2\pi i/p$. If we let $H$ be the subgroup generated by $B$ and elements of $\mathrm{S}^1$ we may define*

$$\chi_i = \begin{cases} \mathrm{Cor}_H^{\widetilde{P}}(\tau'^i) & \text{for } i < p-1 \\ \mathrm{Cor}_H^{\widetilde{P}}(\tau'^{p-1}) - \alpha^{p-1} & \text{for } i = p-1 \end{cases}$$

*where $\tau'$ is any element of $H^2(\mathrm{B}H;\mathbb{Z})$ restricting to $\mathrm{S}^1$ as the standard generator $\tau$. Similarly, $\zeta$ may be defined to be the $p$th Chern class of a $p$-dimensional irreducible representation of $\widetilde{P}$ that restricts to $\mathrm{S}^1$ as $p$ copies of the identity.* ∎

**Theorem 2.** *$H^*(\mathrm{B}\widetilde{P};\mathbb{F}_p)$ is generated by elements $y, y', x, x', c_2, \ldots, c_{p-1}, z$ with*

$$\deg(y) = \deg(y') = 1, \quad \deg(x) = \deg(x') = 2, \quad \deg(c_i) = 2i, \quad \deg(z) = 2p,$$

*subject to the following relations:*

$$xy' = x'y,$$

$$x^p y' = x'^p y, \qquad x^p x' = x'^p x,$$
$$\beta(y) = x, \qquad \beta(y') = x',$$

$$c_i y = \begin{cases} 0 \\ -x^{p-1} y \end{cases} \qquad c_i y' = \begin{cases} 0 & \text{for } i < p-1 \\ -x'^{p-1} y' & \text{for } i = p-1, \end{cases}$$

$$c_i x = \begin{cases} 0 \\ -x^p \end{cases} \qquad c_i x' = \begin{cases} 0 & \text{for } i < p-1 \\ -x'^p & \text{for } i = p-1, \end{cases}$$

$$c_i c_j = \begin{cases} 0 & \text{for } i+j < 2p-2 \\ x^{2p-2} + x'^{2p-2} - x^{p-1} x'^{p-1} & \text{for } i = j = p-1. \end{cases}$$

If we let $\pi_*$ stand for the projection map from $H^n(\mathrm{B}G; \mathbb{Z})$ to $H^n(\mathrm{B}G; \mathbb{F}_p)$, and $\delta_p$ for the Bockstein from $H^n(\mathrm{B}G; \mathbb{F}_p)$ to $H^{n+1}(\mathrm{B}G; \mathbb{Z})$, then we may define the generators by the following equations.

$$\pi_*(\chi_i) = c_i, \quad \pi_*(\zeta) = z, \quad \pi_*(\alpha) = x, \quad \pi_*(\beta) = x',$$

$$\delta_p(y) = \alpha, \quad \delta_p(y') = \beta.$$

This determines the effect of automorphisms of $\widetilde{P}$ on the cohomology ring. The element $yy'$ is a non-zero multiple of $\pi_*(\chi_1)$.

*Proof.* From Theorem 1 we see that $H^n(\mathrm{B}\widetilde{P}; \mathbb{Z})$ is trivial for $n$ odd, and has torsion of exponent at most $p$ for $n$ even. It follows that $\pi_*$ maps $H^{2n}(\mathrm{B}\widetilde{P}; \mathbb{Z})$ onto $H^{2n}(\mathrm{B}\widetilde{P}; \mathbb{F}_p)$, with kernel the elements divisible by $p$, and that $\delta_p$ maps $H^{2n-1}(\mathrm{B}\widetilde{P}; \mathbb{F}_p)$ isomorphically to the torsion subgroup of $H^{2n}(\mathrm{B}\widetilde{P}; \mathbb{Z})$. This determines the ring structure of $H^{\mathrm{even}}(\mathrm{B}\widetilde{P}; \mathbb{F}_p)$, and the module structure for this ring of $H^{\mathrm{odd}}(\mathrm{B}\widetilde{P}; \mathbb{F}_p)$. More explicitly, the relations of odd degree given in the statement must hold, because they map under the injective map $\delta_p$ to valid relations, for example $\delta_p(xy') = \alpha\beta = \delta_p(x'y)$. No more generators of odd degree are required, because the torsion in $H^{\mathrm{even}}(\mathrm{B}\widetilde{P}; \mathbb{Z})$ is spanned by multiples of $\alpha$ and $\beta$. The given relations of odd degree suffice to show that monomials of the forms $z^i x^j y$ and $z^i x^j x'^k y'$, where either $k \leqslant p-2$ or $j = 0$, span $H^{\mathrm{odd}}(\mathrm{B}\widetilde{P}; \mathbb{F}_p)$, and $\delta_p$ maps these to a basis for the torsion in $H^*(\mathrm{B}\widetilde{P}; \mathbb{Z})$, so no more relations of odd degree are needed. Since $yy'x = yy'x' = 0$, it follows that $yy'$ must be some multiple of $\pi_*(\chi_1)$, but it remains to show that $yy'$ is non-zero. The easiest way to check this is by examining the spectral sequence for $\widetilde{P}$ expressed as an extension of $\mathrm{S}^1$ by $C_p \oplus C_p$. In this spectral sequence $y$ and $y'$ yield generators for $E_2^{1,0}$, their product is non-zero in $E_2^{2,0}$, and $E_2^{0,1}$ is trivial, so $E_2^{2,0} = E_\infty^{2,0}$. Alternatively, it may be verified using explicit cocycles that $\pi_*(\chi_1) = y'y$, which must therefore be non-zero. ∎

Since $H^*(\mathrm{B}\widetilde{P}; \mathbb{F}_p)$ is generated by elements of degree at most $2p$, and we already know the action of the Bockstein, the following proposition suffices to determine the action of the Steenrod algebra on $H^*(\mathrm{B}\widetilde{P}; \mathbb{F}_p)$.

Proposition 3. *With notation as in Theorem 2, and writing $c_1$ for $y'y$, the following relations hold.*

$$\mathrm{P}^1(z) = zc_{p-1}, \qquad \mathrm{P}^1(c_i) = \begin{cases} izc_{i-1} & \text{for } i < p-1, \\ -zc_{p-2} + x^{2p-2} + x'^{2p-2} - x^{p-1}x'^{p-1} & \text{for } i = p-1. \end{cases}$$

*Proof.* As in the statement of Theorem 1, let $H$ be the subgroup of $\widetilde{P}$ generated by $B$ and $S^1$. Let $t$ be the image under $\pi_*$ of the element $\tau'$ of $H^2(BH;\mathbb{Z})$, write Cor and Res for the corestriction and restriction maps between $\widetilde{P}$ and $H$, and let $\bar{x}' = \mathrm{Res}(x')$. Then the definition of $c_i$ given above is equivalent to $c_i = \mathrm{Cor}(t^i)$ for $i < p-1$, and $c_{p-1} = \mathrm{Cor}(t^{p-1}) - x^{p-1}$. The $p$-dimensional representation of $\widetilde{P}$ with $p$th Chern class $z$ restricts to $H$ as a sum of one copy of each of the representations of $H$ that restrict to $S^1$ as the identity. The first Chern classes of these representations are $t + j\bar{x}'$ for each choice of $j$, and so

$$\mathrm{Res}(z) = \prod_{j=0}^{p-1}(t + j\bar{x}') = t^p - \bar{x}'^{p-1}t.$$

The action of conjugation by $A$ on $H^1(BH;\mathbb{F}_p)$ sends $t$ to $t+\bar{x}'$ and fixes $\bar{x}'$. The restriction-corestriction formula [3] therefore implies that

$$\mathrm{ResCor}(t^i) = \sum_{j=0}^{p-1}(t+j\bar{x}')^i = \sum_{k=0}^{i}\binom{i}{k}t^{i-k}\bar{x}'^k\left(\sum_{j=0}^{p-1}j^k\right).$$

Newton's formula states that $\sum_{j=1}^{p-1}j^k$ is divisible by $p$ except when $k$ is divisible by $p-1$, in which case it is congruent to $-1$ mod $p$. Applying this formula we see that

$$\mathrm{ResCor}(t^i) = \begin{cases} 0 & \text{for } i < p-1 \\ -\bar{x}'^{p-1} & \text{for } i = p-1. \end{cases}$$

Now we may determine $\mathrm{P}^1(c_i)$ by repeated application of Frobenius reciprocity, as below.

$$\begin{aligned}\mathrm{P}^1\mathrm{Cor}(t^i) &= \mathrm{Cor}(\mathrm{P}^1(t^i)) = i\mathrm{Cor}(t^{p+i-1}) \\ &= i\mathrm{Cor}(t^{i-1}\mathrm{Res}(z) + t^i\bar{x}'^{p-1}) \\ &= iz\mathrm{Cor}(t^{i-1}) + ix^{p-1}\mathrm{Cor}(t^i).\end{aligned}$$

To check that $\mathrm{P}^1(z)$ is as claimed, note that a general element of $H^{4p-2}(B\widetilde{P};\mathbb{F}_p)$ has the form $\lambda zc_{p-1} + zP(x,x') + Q(x,x')$, so we may assume that $\mathrm{P}^1(z)$ has that form. First we show that $Q(x,x')$ is zero. Let $H'$ be any subgroup of $\widetilde{P}$ of index $p$, and let $\hat{t}, \hat{x} \in H^2(BH;\mathbb{F}_p)'$ be such that $\hat{t}$ restricts to a generator for $H^2(BS^1;\mathbb{F}_p)$ and $\hat{x}$ restricts trivially. As in the case $H' = H$ considered above, it may be shown that $\mathrm{Res}_{H'}^{\widetilde{P}}(z)$, and hence also $\mathrm{Res}_{H'}^{\widetilde{P}}(\mathrm{P}^1z)$, are divisible by $\hat{t}$, whereas $Q(x,x')$ restricts to $H'$ as a scalar multiple of $\hat{x}^{2p-1}$. This multiple is zero for each choice of $H'$ if and only if $Q(x,x')$ is a multiple of each $x + ix'$ and of $x'$. In this case $Q$ is a multiple of $x^p x' - x'^p x$ so represents zero in $H^{4p-2}(B\widetilde{P}; \mathbb{F}_p)$.

To show that $P(x, x')$ is zero, note that $z$ and $\mathrm{P}^1(z)$ are invariant under the action of $SL_2(\mathbb{F}_p)$ considered as a subgroup of $\mathrm{Aut}(\widetilde{P})$, but no non-zero $P(x, x')$ of degree $p-1$ is invariant. Thus we have shown that $\mathrm{P}^1(z)$ is a multiple of $zc_{p-1}$. This multiple may be determined by restricting to the subgroup $H$. ∎

We can now state our result for $H^*(\mathrm{B}P(n); \mathbb{F}_p)$ for $n \geqslant 4$.

**Theorem 4.** *Let $P(n)$ be the group presented in the introduction, and view it as a subgroup of $\widetilde{P}$. Then for $n \geqslant 4$ the restriction map from $H^*(\mathrm{B}\widetilde{P}; \mathbb{F}_p)$ to $H^*(\mathrm{B}P(n); \mathbb{F}_p)$ is injective, and if we write $\Lambda[u]$ for an exterior algebra over $\mathbb{F}_p$ on a generator $u$ of degree one, then $H^*(\mathrm{B}P(n); \mathbb{F}_p)$ is isomorphic to the graded tensor product $\Lambda[u] \otimes H^*(\mathrm{B}P(n); \mathbb{F}_p)$. We may take for $u$ the element given by the cocycle $u([A^i B^j C^k]) = k$. Since any automorphism of $P(n)$ extends to an endomorphism of $\widetilde{P}$, Theorem 2 determines the action of $\mathrm{Aut}(P(n))$ on $H^*(\mathrm{B}P(n); \mathbb{F}_p)$. The action of the Steenrod algebra is determined by Proposition 3, except for $\beta(u)$, which is zero for $n > 4$ and $y'y$ for $n = 4$.*

*Proof.* The $E_2$ page of the spectral sequence for $\mathrm{B}P(n)$ as an $\mathrm{S}^1$ bundle over $\mathrm{B}\widetilde{P}$ is isomorphic to $\Lambda[u] \otimes H^*(\mathrm{B}P(n); \mathbb{F}_p)$, so $E_2^{1,0} \cong C_p \oplus C_p$ and $E_2^{0,1} \cong C_p$. We know that $H^1(\mathrm{B}P(n); \mathbb{F}_p) \cong \mathrm{Hom}(P(n), C_p) \cong (C_p)^3$, so $u$ must survive, and the spectral sequence collapses. It remains to check that $\beta(u)$ is as claimed. The image of $\delta_p$ from $H^1(\mathrm{B}P(n); \mathbb{F}_p)$ to $H^2(\mathrm{B}P(n); \mathbb{Z})$ is the subgroup of elements of order $p$, and the kernel of $\pi_*$ from $H^2(\mathrm{B}P(n); \mathbb{Z})$ to $H^2(\mathrm{B}P(n); \mathbb{F}_p)$ is the subgroup of elements divisible by $p$. Since $H^2(\mathrm{B}P(n); \mathbb{Z}) \cong \mathrm{Hom}(P(n), \mathrm{S}^1) \cong C_{p^{n-3}} \oplus C_p \oplus C_p$, it follows that $u$ may be chosen such that $\beta(u)$ is as claimed. It may be checked using cochains that $u$ as described in the statement is such an element. ∎

In the case when $n = 3$ the spectral sequence considered above does not collapse, and many generators not in the image of the restriction from $\widetilde{P}$ (or equivalently, yielding elements on the $E_\infty^{*,1}$ line of the spectral sequence) are required. Some of these may be described as corestrictions from subgroups of $P(3)$, but others require us to consider Massey products. Many of the relations between these elements will be implied by relations that Massey products always satisfy, which we summarise in the following lemma.

**Lemma 5.** *Let $C^*$ stand for the cochain complex for a space $X$ with $\mathbb{F}_p$ coefficients, and let $H^*$ be the cohomology of $C^*$. Given homogeneous cocycles $u, v, w$ in $C^*$ such that $[uv]$ and $[vw]$ represent 0 in $H^*$, choose cochains $a$ and $b$ such that $\delta a = uv$, $\delta b = vw$, and define the Massey product of $[u]$, $[v]$ and $[w]$ by the formula*

$$\langle [u], [v], [w] \rangle = [(-1)^{|u|} ub - aw] \in H^{|u|+|v|+|w|-1} / (uH^{|v|+|w|-1} + wH^{|u|+|v|-1}).$$

*This product is linear in each of its arguments, and satisfies the following relations whenever all the terms are defined.*

$$\langle u, v, w \rangle x + (-1)^{|u|} u \langle v, w, x \rangle \equiv 0 \quad \mathrm{mod}\ uH^* x \qquad (1)$$

$$(-1)^{|u|} \langle \langle u, v, w \rangle, x, y \rangle + \langle u, \langle v, w, x \rangle, y \rangle + (-1)^{|v|} \langle u, v, \langle w, x, y \rangle \rangle \equiv 0$$
$$\mathrm{mod}\ uH^* + H^{|u|+|v|-1} w H^{|x|+|y|-1} + yH^* \qquad (2)$$

$$(-1)^{|w||u|}\langle u,v,w\rangle + (-1)^{|u||v|}\langle v,w,u\rangle + (-1)^{|v||w|}\langle w,u,v\rangle \equiv 0$$
$$\mod uH^* + vH^* + wH^* \quad (3)$$

$$\langle u,v,w\rangle + (-1)^{|u||v|+|v||w|+|w||u|}\langle w,v,u\rangle \equiv 0 \mod uH^* + wH^* \quad (4)$$

*Proof.* See for example [8] or [12]. ∎

We are ready now to calculate $H^*(\mathrm{B}P(3); \mathbb{F}_p)$, which we split into the cases $p > 3$ and $p = 3$. It should not be surprising that these cases are quite different, because the groups of order $p^4$ behave differently in these cases. For example, there are two generator groups of order $p^4$ and exponent $p$ if and only if $p > 3$ [4]. Any such group is expressible as a non-split central extension of $C_p$ by $P(3)$, such that every 'subextension' of $C_p$ by $C_p$ is split. Translating these statements into a statement about extension classes in $H^2(\mathrm{B}P(3); \mathbb{F}_p)$, we see that $H^2(\mathrm{B}P(3); \mathbb{F}_p)$ is detected by cyclic subgroups if and only if $p = 3$.

**Theorem 6** *Let $p$ be greater than 3. Then $H^*(\mathrm{B}P(3); \mathbb{F}_p)$ is generated by elements $y, y'$, $x, x', Y, Y', X, X', d_4, \ldots, d_p, c_4, \ldots, c_{p-1}, z$, with*

$$\deg(y) = \deg(y') = 1, \quad \deg(x) = \deg(x') = \deg(Y) = \deg(Y') = 2,$$

$$\deg(X) = \deg(X') = 3, \quad \deg(d_i) = 2i - 1, \quad \deg(c_i) = 2i, \quad \deg(z) = 2p,$$

$$\beta(y) = x, \quad \beta(y') = x',$$
$$\beta(Y) = X, \quad \beta(Y') = X',$$

$$\beta(d_i) = \begin{cases} c_i & \text{for } i < p \\ 0 & \text{for } i = p, \end{cases}$$

*subject to the following relations:*

$$yy' = 0, \quad xy' = x'y,$$

$$yY = y'Y' = 0, \quad yY' = y'Y,$$

$$Y^2 = Y'^2 = YY' = 0,$$

$$yX = xY, \quad y'X' = x'Y',$$
$$Xy' = 2xY' + x'Y, \quad X'y = 2x'Y + xY',$$

$$XY = X'Y' = 0, \quad XY' = -X'Y, \quad xX' = -x'X,$$

$$x(xY' + x'Y) = x'(xY' + x'Y) = 0,$$

$$x^p y' - x'^p y = 0, \quad x^p x' - x'^p x = 0,$$
$$x^p Y' + x'^p Y = 0, \quad x^p X' + x'^p X = 0,$$

$$c_i y = \begin{cases} 0 \\ -x^{p-1}y \end{cases} \qquad c_i y' = \begin{cases} 0 & \text{for } i < p-1 \\ -x'^{p-1}y' & \text{for } i = p-1, \end{cases}$$

$$c_i x = \begin{cases} 0 \\ -x^p \end{cases} \qquad c_i x' = \begin{cases} 0 & \text{for } i < p-1 \\ -x'^p & \text{for } i = p-1, \end{cases}$$

$$c_i Y = \begin{cases} 0 \\ -x^{p-1}Y \end{cases} \qquad c_i Y' = \begin{cases} 0 & \text{for } i < p-1 \\ -x'^{p-1}Y' & \text{for } i = p-1, \end{cases}$$

$$c_i X = \begin{cases} 0 \\ -x^{p-1}X \end{cases} \qquad c_i X' = \begin{cases} 0 & \text{for } i < p-1 \\ -x'^{p-1}X' & \text{for } i = p-1, \end{cases}$$

$$c_i c_j = \begin{cases} 0 & \text{for } i+j < 2p-2 \\ x^{2p-2} + x'^{2p-2} - x^{p-1}x'^{p-1} & \text{for } i = j = p-1, \end{cases}$$

$$d_i y = \begin{cases} 0 \\ -x^{p-1}Y \end{cases} \qquad d_i y' = \begin{cases} 0 & \text{for } i < p \\ x'^{p-1}Y' & \text{for } i = p, \end{cases}$$

$$d_i x = \begin{cases} 0 \\ -x^{p-1}y \\ x^{p-1}X \end{cases} \qquad d_i x' = \begin{cases} 0 & \text{for } i < p-1 \\ -x'^{p-1}y' & \text{for } i = p-1 \\ -x'^{p-1}X' & \text{for } i = p, \end{cases}$$

$$d_i Y = 0 \qquad d_i Y' = 0,$$

$$d_i X = \begin{cases} 0 \\ -x^{p-1}Y \end{cases} \qquad d_i X' = \begin{cases} 0 & \text{for } i \neq p-1 \\ -x'^{p-1}Y' & \text{for } i = p-1, \end{cases}$$

$$d_i d_j = \begin{cases} 0 & \text{for } i < p-1 \text{ or } j < p-1 \\ x^{2p-3}Y - x'^{2p-3}Y' + x^{p-1}x'^{p-2}Y' & \text{for } i = p \text{ and } j = p-1, \end{cases}$$

$$d_i c_j = \begin{cases} 0 & \text{for } i < p-1 \text{ or } j < p-1 \\ x^{2p-3}y + x'^{2p-3}y' - x^{p-1}x'^{p-2}y' & \text{for } i = j = p-1 \\ -x^{2p-3}X + x'^{2p-3}X' - x^{p-1}x'^{p-2}X' & \text{for } i = p, j = p-1. \end{cases}$$

We define $y$, $y'$ by the cocycles $y(A^i B^j C^k) = i$ and $y'(A^i B^j C^k) = j$. The equation $yy' = 0$ implies that we may define unique elements $Y$, $Y'$ by

$$Y = \langle y, y, y' \rangle, \qquad Y' = \langle y', y', y \rangle.$$

If we define $d' \in H^1(\mathrm{B}\langle B, C \rangle; \mathbb{F}_p)$ by the cocycle $d'(B^r C^s) = s$, and $c' = \beta(d')$, then we may define

$$d_i = \begin{cases} \mathrm{Cor}_{\langle B,C \rangle}^{P(3)}(c'^{i-1}d') & \text{for } i < p-1 \\ \mathrm{Cor}_{\langle B,C \rangle}^{P(3)}(c'^{p-2}d') - x^{p-2}y & \text{for } i = p-1 \\ \mathrm{Cor}_{\langle B,C \rangle}^{P(3)}(c'^{p-1}d') + x^{p-2}X & \text{for } i = p. \end{cases}$$

We may define the elements $x$, $x'$, $X$, $X'$, $c_i$ from the elements already defined, using the Bockstein, and we may define the elements $y$, $y'$, $x$, $x'$, $c_i$, $z$ to be the restrictions from $\widetilde{P}$ of the generators with the same name. The two definitions given for $y$, $y'$, $x$, $x'$, $c_i$ are equivalent. The effect of automorphisms of $P(3)$ on $y$, $y'$, $x$, $x'$, $Y$, $Y'$, $X$, $X'$ is determined by the above definitions. An automorphism of $P(3)$ that restricts to its centre as $C \mapsto C^j$ sends $d_i$ to $j^i d_i$, $c_i$ to $j^i c_i$ and $z$ to $jz$. We note that the restriction from $\widetilde{P}$ to $P(3)$ sends $c_2$ to a non-zero multiple of $xY' + x'Y$, and $c_3$ to a non-zero multiple of $XX'$.

*Proof.* Throughout this proof we shall adopt the convention that Cor and Res will stand for $\mathrm{Cor}^{P(3)}_{\langle B,C \rangle}$ and $\mathrm{Res}^{P(3)}_{\langle B,C \rangle}$ respectively. As in the previous theorem we consider the spectral sequence for $P(3)$ as an $S^1$ bundle over $\widetilde{P}$. First we prove that $H^n(BP(3); \mathbb{Z})$ has exponent $p$ for $n < 2p$ by considering the spectral sequence with integer coefficients. (This is proved in [11] and [10], but we include it for completeness.) Since $H^{2k+1}(B\widetilde{P}; \mathbb{Z})$ is trivial, it follows that $H^{2k+1}(BP(3); \mathbb{Z})$ is isomorphic to a subgroup of $H^{2k}(B\widetilde{P}; \mathbb{Z})$, which has torsion only of order $p$. It also follows that $H^{2k}(B\widetilde{P}; \mathbb{Z})$ maps onto $H^{2k}(BP(3); \mathbb{Z})$. Since for $k < p$ $H^{2k}(B\widetilde{P}; \mathbb{Z})$ is generated by elements of order $p$ and corestrictions from proper subgroups the same must be true for $H^{2k}(BP(3); \mathbb{Z})$, which therefore has exponent $p$. It follows that for $n < 2p$ the image and kernel of the Bockstein $\beta$ in $H^n(BP(3); \mathbb{F}_p)$ coincide.

Returning now to the spectral sequence with $\mathbb{F}_p$ coefficients, the $E_2$ page is isomorphic to $\Lambda[u] \otimes H^*(B\widetilde{P}; \mathbb{F}_p)$, where $u$ has bidegree $(0, 1)$. Since $H^1(BP(3); \mathbb{F}_p)$ has order $p^2$, $d_2(u)$ must be non-zero. All the subgroups of $P(3)$ of index $p$ are isomorphic to $C_p \oplus C_p$, so $d_2(u)$ must restrict trivially to each $S^1 \oplus C_p$ subgroup of $\widetilde{P}$. It follows that $d_2(u)$ is a non-zero multiple of $yy'$, and so

$$E^{*,0}_\infty \cong H^*(B\widetilde{P}; \mathbb{F}_p)/(yy'), \qquad E^{*,1}_\infty \cong \mathrm{Ann}(yy') \subset H^*(B\widetilde{P}; \mathbb{F}_p).$$

Hence $E^{*,*}_\infty$ is generated by the elements $[uy], [uy'], [ux], [ux'], [uc_2], \ldots, [uc_{p-1}], y, y', x, x'$, $c_2, \ldots, c_{p-1}, z$, subject to the relations implied by those that hold in $E^{*,*}_2$, together with the relation $yy' = 0$. The cohomology relations that involve only generators in the image of the restriction from $\widetilde{P}$ follow from the spectral sequence.

Since unique elements $Y$, $Y'$ may be defined as in the statement, our next task is to show that these elements together with $x$, $x'$ form a basis for $H^2$. We may calculate $\mathrm{Res}^{P(3)}_{\langle A,C \rangle}(Y')$ and $\mathrm{Res}(Y') = \mathrm{Res}^{P(3)}_{\langle B,C \rangle}(Y')$ by explicit calculation using cochains in the bar resolution for $P(3)$. Define 1-cochains $a$ and $b$ by the equations $a([g]) = -1/2(y'([g]))^2$, $b([A^r B^s C^t]) = t$. Then

$$\delta a([g|h]) = -1/2 \left( y'([g])^2 + y'([h])^2 - (y'([g]) + y'([h]))^2 \right) = y'([g])y'([h]), \text{ and}$$

$$\delta b([A^i B^j C^k | A^r B^s C^t]) = b([A^i B^j C^k]) + b([A^r B^s C^t]) - b([A^{i+r} B^{j+s} C^{k+t-jr}])$$
$$= y'([A^i B^j C^k])y([A^r B^s C^t]),$$

so we may choose a cocycle representing $Y'$ by the formula

$$Y'([g|h]) = y'([g])(1/2 y'([g])y([h]) - b([h])).$$

The cocycle $y$ vanishes on $\langle B, C \rangle$, and $y'$ vanishes on $\langle A, C \rangle$, so it is now easy to check that even as cochains $\mathrm{Res}^{P(3)}_{\langle A,C \rangle}(Y') = 0$ and $\mathrm{Res}(Y') = \mathrm{Res}(y')d'$, where $d'$ is the cocycle defined in the statement. A similar argument may be used to show that $Y$ restricts to zero on $\langle B, C \rangle$ and to a non-zero product of elements of degree one on $\langle A, C \rangle$. Since $x$ and $x'$ restrict to subgroups as elements in the image of the Bockstein, it follows that $x$, $x'$, $Y$, $Y'$ form a basis for $H^2$. The identities of Lemma 5 suffice to express any Massey product of elements of $H^1$ in terms of $Y$ and $Y'$, so this determines the action of $\mathrm{Aut}(P(3))$ on $Y$ and $Y'$. For example, the 'Jacobi identity' ((3) of Lemma 5) implies that $3\langle y, y, y\rangle = 0$.

The relations in $H^3$ between $Y$, $Y'$ and $y$, $y'$ follow from properties of the Massey product, for example

$$Yy' = \langle y, y, y'\rangle y' \equiv y\langle y, y', y'\rangle \equiv y\langle y', y', y\rangle = yY',$$

where the congruences are modulo $\{0\}$, and come from applications of (1) and (4) of Lemma 5. Before moving on to $H^4$ we note that $\langle y, y, Y\rangle$ is defined, and is congruent to zero modulo $H^2 H^1$, which follows from (2) of Lemma 5. Applying each of (1) and (4) it is now easy to show that $Y^2 \equiv 0$ modulo $H^2 H^1 H^1 = \{0\}$. Similarly $Y'^2 = (Y+Y')^2 = 0$, because there are automorphisms of $P(3)$ sending $Y$ to $Y'$ and to $Y+Y'$. The expressions for $yX$ and $y'X'$ follow by applying the Bockstein to the relations $yY = y'Y' = 0$. For the remaining relations in $H^4$ we introduce some matrix Massey products [12]. We consider $\langle (x, x'), \begin{pmatrix} y' \\ -y \end{pmatrix}, y\rangle$, which is defined modulo $H^2 H^1$, so has a well defined product with any element of $H^1$. These products are determined by the following calculations.

$$\langle (x, x'), \begin{pmatrix} y' \\ -y \end{pmatrix}, y\rangle y \equiv -(x, x')\begin{pmatrix} \langle y', y, y\rangle \\ -\langle y, y, y\rangle \end{pmatrix} = -(x, x')\begin{pmatrix} Y \\ 0 \end{pmatrix} = -xY,$$

$$\langle (x, x'), \begin{pmatrix} y' \\ -y \end{pmatrix}, y\rangle y' \equiv -(x, x')\begin{pmatrix} \langle y', y, y'\rangle \\ -\langle y, y, y'\rangle \end{pmatrix} \equiv -(x, x')\begin{pmatrix} -2Y' \\ -Y \end{pmatrix} = 2xY' + x'Y,$$

where each congruence is modulo $H^2 H^1 H^1 = \{0\}$. Similarly, the following expressions may be verified.

$$\langle (x', x), \begin{pmatrix} y \\ -y' \end{pmatrix}, y'\rangle y' = -x'Y', \qquad \langle (x', x), \begin{pmatrix} y \\ -y' \end{pmatrix}, y'\rangle y = 2x'Y + xY'.$$

Hence we deduce that $X \in \langle (x, x'), \begin{pmatrix} y' \\ -y \end{pmatrix}, y\rangle$ and $X' \in \langle (x', x), \begin{pmatrix} y \\ -y' \end{pmatrix}, y'\rangle$, and the remaining relations in $H^4$ follow. Our results allow us to deduce that $X, X', yY', xy, xy', x'y'$ form a basis for $H^3$. In $H^5$ the relations $XY = 0$, $X'Y' = 0$, $XY' = -X'Y$ and $xX' = -x'X$ follow by applying the Bockstein to the relations $Y^2 = 0$, $Y'^2 = 0$, $YY' = 0$ and $X'y = 2x'Y + xY'$ respectively. For the relations stated in $H^6$, we note that $xyY' = x'yY' = 0$, and then apply the Bockstein to these relations, noting also that $\beta(yY') = 2xY' + 2x'Y$.

We see now that in the $E^{*,*}_\infty$ page of the spectral sequence $Y$ yields $\lambda[uy]$ and $Y'$ yields $-\lambda[uy']$ for some non-zero $\lambda$, and deduce that $xY' + x'Y \equiv 0$ modulo the image of the restriction from $\widetilde{P}$. The relations in $H^6$ imply that $xY' + x'Y = \lambda' \mathrm{Res}^{\widetilde{P}}_{P(3)}(c_2)$ for some $\lambda'$. The element $yY'$ is not in the kernel of the Bockstein, so $\lambda'$ must be non-zero, and we see that products of elements in $H^2$ generate $H^4$.

The spectral sequence shows that $xH^3 + x'H^3$ has index $p$ in $H^5$, and it may be shown that $XY' \notin xH^3 + x'H^3$ by considering another spectral sequence, that for $P(3)$ expressed as a central extension of $C_p$ by $C_p \oplus C_p$. Only the first two differentials need be computed for this purpose. This calculation is contained in [10]. Since $XY'$ is therefore not in the image of the Bockstein it follows that $\beta(XY') = -XX'$ is not zero. $XX'$ is however annihilated by $x$, so must be a multiple of $c_3$. Hence $H^6$ is generated by $XX'$ and products of elements of $H^2$.

We already have the relations $x^p y' = x'^p y$, $x^p x' = x'^p x$, and we obtain $x^p Y' + x'^p Y = 0$ by applying $\mathrm{P}^1$ to the relation $xY' + x'Y = \lambda \mathrm{Res}_{P(3)}^{\widetilde{P}}(c_2)$. (Recall that $\mathrm{P}^1(c_2)$ was determined in Proposition 3.) The relation $x^p X' + x'^p X = 0$ now follows by applying the Bockstein.

We may verify that in degrees greater than 6 all products of the generators $y, y', x, x'$, $Y, Y', X, X'$ (which we shall call the 'low dimensional generators') may be expressed in the form
$$f_1 + f_2 Y + f_3 Y' \qquad \text{for even total degree}$$
$$f_1 y + f_2 y' + f_3 X + f_4 X' \qquad \text{for odd total degree}$$
where $f_i$ is a polynomial in $x$ and $x'$. With the exception of $xY' + x'Y$, such expressions satisfy 'the same' relations as elements of $H^*$ as they do as elements of $E_\infty^{*,*}$. Elements that are expressible as above form a subspace of $H^n$ of codimension 1 for $7 \leqslant n \leqslant 2p$, so we introduce the elements $c_i$ for $i > 3$ and $z$ to our generating set. The Bockstein $\beta$ sends elements of the above form of odd total degree to elements of the above form of even total degree, so any element sent by $\beta$ to $c_i$ will suffice to complete a basis for $H^{2i-1}$. Hence we may add $d_i$ for $i < p$ to our generating set. Assuming that the relations given involving $d_p$ hold, it follows that $d_p$ is not an element of the form described above, so suffices to complete a basis for $H^{2p-1}$. Using Frobenius reciprocity and the fact that CorRes is the zero map, it may be checked that $\beta(d_p) = 0$.

The relations involving the $c_i$ and $d_i$ may be checked using Frobenius reciprocity, and we shall only prove some examples and indicate an economical order in which to prove the rest. Before starting it is helpful to notice that the automorphism $\theta$ of $P(3)$ defined by $\theta(A) = B$, $\theta(B) = A$ has the effect of exchanging the 'primed' and 'unprimed' low dimensional generators, and sends $c_i$ to $(-1)^i c_i$.

The relations between $c_i$ and $y, y', x, x'$ follow from the spectral sequence. The expressions for $c_i Y$ are easily determined since $\mathrm{Cor}(c'^i)Y = \mathrm{Cor}(c'^i \mathrm{Res}(Y)) = 0$. The expressions for $c_i Y'$ follow by applying $\theta^*$, and the expressions for $c_i X$ and $c_i X'$ by applying $\beta$. The relations between $d_i$ and the low dimensional generators are slightly more involved, and we must examine the 'primed' and 'unprimed' generators separately, because we do not know yet what effect $\theta^*$ has on $d_i$. For example, the following calculations verify the claimed expressions for $d_i y$ and $d_i y'$.
$$\mathrm{Cor}(c'^{i-1} d') y = \mathrm{Cor}(c'^{i-1} d' \mathrm{Res}(y)) = 0,$$
$$\mathrm{Cor}(c'^{i-1} d') y' = \mathrm{Cor}(c'^{i-1} d' \mathrm{Res}(y')) = -\mathrm{Cor}(c'^{i-1}) Y'.$$
Similar calculations may be used to determine $d_i Y$ and $d_i Y'$, and then we may apply the Bockstein to these relations (substituting from earlier relations for terms involving $c_i$) to obtain expressions for $d_i x$, $d_i x'$, $d_i X$ and $d_i X'$.

The expressions for $d_i d_j$ and $c_i d_j$ may be determined similarly. As an example we determine $\mathrm{Cor}(c'^{i-1}d')\mathrm{Cor}(c'^{j-1}d')$, where without loss of generality we may assume that $i > j$. If we write $\bar{y}' = \mathrm{Res}(y')$, $\bar{x}' = \mathrm{Res}(x')$, then the action of conjugation by $A$ on $H^1(\mathrm{B}\langle B, C\rangle; \mathbb{F}_p)$ fixes $\bar{y}'$ and sends $d'$ to $d' + \bar{y}'$. It follows from the restriction-corestriction formula that

$$\mathrm{ResCor}(c'^{j-1}d') = \sum_{k=0}^{p-1} (c' + k\bar{x}')^{j-1}(d' + k\bar{y}').$$

Applying Newton's formula as in the proof of Proposition 3 we see that this expression is zero for $j < p-1$, and that $\mathrm{ResCor}(c'^{p-2}d') = -\bar{y}'\bar{x}'^{p-2}$. It only remains to determine $\mathrm{Cor}(c'^{p-1}d')\mathrm{Cor}(c'^{p-2}d')$, which may be determined as below.

$$\begin{aligned}
\mathrm{Cor}(c'^{p-1}d')\mathrm{Cor}(c'^{p-2}d') &= -\mathrm{Cor}(c'^{p-1}d'\bar{y}'\bar{x}'^{p-2}) \\
&= \mathrm{Cor}(c'^{p-1})Y'\bar{x}'^{p-2} \\
&= (c_{p-1} + x^{p-1})Y'x'^{p-2} \\
&= x^{p-1}x'^{p-2}Y' - x'^{2p-3}Y'.
\end{aligned}$$

All that remains to be calculated is the effect of automorphisms of $P(3)$ on $d_i$. The multiplicative relations demonstrated above imply that in the spectral sequence $d_i$ yields the element $[uc_{i-1}]$ in $E_\infty^{2i-2,1}$. An automorphism of $P(3)$ that restricts to the centre as the map $C \mapsto C^j$ extends to an endomorphism of $\widetilde{P}$ that wraps the $S^1$ subgroup $j$ times around itself. This endomorphism induces an endomorphism of the spectral sequence which sends $u$ to $ju$, and $c_i$ to $j^i c_i$. It follows that $d_i$ is sent to $j^i d_i$ modulo $E_\infty^{2i-1,0}$. Since $\beta(d_i) = c_i$ for $i < p$ and $\beta(d_p) = 0$, we also know that $d_i$ is sent to $j^i d_i$ modulo the kernel of $\beta$. The result now follows, since $\mathrm{Ker}(\beta) \cap E_\infty^{2i-1,0}$ is trivial. ∎

We now state and prove a similar assertion for the case when $p = 3$. This result has also been obtained by Milgram and Tezuka [13], who show that the whole ring is detected by the restriction maps to proper subgroups. For $p > 3$ there are elements (for example $XX'$, $c_i$ for $i < p-1$) that restrict trivially to all proper subgroups of $P(3)$. One way to explain the difference between the presentations of $H^*(\mathrm{B}P(3); \mathbb{F}_p)$ in the cases $p > 3$ and $p = 3$ is in terms of the different properties of Massey products. If $A$ is any space, and $v \in H^1(A; \mathbb{F}_p)$, then $\langle v, v, v \rangle = 0$ if $p > 3$, whereas $\langle v, v, v \rangle = \beta(v)$ if $p = 3$. This may be verified using explicit cochain calculations for the 'universal example' which is the case when $A = \mathrm{B}C_p$, or see [8].

**Theorem 7.** Let $p = 3$. Then $H^*(\mathrm{B}P(3); \mathbb{F}_3)$ is generated by elements $y, y', x, x', Y, Y', X, X', z$, with

$$\deg(y) = \deg(y') = 1, \quad \deg(x) = \deg(x') = \deg(Y) = \deg(Y') = 2,$$

$$\deg(X) = \deg(X') = 3, \quad \deg(z) = 6,$$

$$\beta(y) = x, \quad \beta(y') = x',$$
$$\beta(Y) = X, \quad \beta(Y') = X',$$

*subject to the following relations:*

$$yy' = 0, \qquad xy' = x'y,$$

$$yY = y'Y' = xy', \qquad yY' = y'Y,$$

$$YY' = xx', \qquad Y^2 = xY', \qquad Y'^2 = x'Y,$$

$$yX = xY - xx', \qquad y'X' = x'Y' - xx',$$
$$Xy' = x'Y - xY', \qquad X'y = xY' - x'Y,$$
$$XY = x'X, \qquad X'Y' = xX',$$
$$XY' = -X'Y, \qquad xX' = -x'X,$$

$$XX' = 0, \qquad x(xY' + x'Y) = -xx'^2, \qquad x'(xY' + x'Y) = -x'x^2,$$

$$x^3 y' - x'^3 y = 0, \qquad x^3 x' - x'^3 x = 0,$$
$$x^3 Y' + x'^3 Y = -x^2 x'^2, \qquad x^3 X' + x'^3 X = 0,$$

We define $y, y' \in H^1(\mathrm{B}P(3); \mathbb{F}_3)$ by the cocycles $y(A^r B^s C^t) = r$ and $y'(A^r B^s C^t) = s$. The equation $yy' = 0$ implies that we may define unique elements in $H^2(\mathrm{B}P(3); \mathbb{F}_3)$ by forming the Massey product of any three elements of $H^1(\mathrm{B}P(3); \mathbb{F}_3)$, and we define $Y, Y'$ by $Y = \langle y, y, y' \rangle$, $Y' = \langle y', y', y \rangle$. We also define $x = \beta(y)$, $x' = \beta(y')$, $X = \beta(Y)$, $X' = \beta(Y')$, and note that $x, x'$ satisfy $x = \langle y, y, y \rangle$, $x' = \langle y', y', y' \rangle$. The effect of automorphisms of $P(3)$ on the generators $y, y', x, x', Y, Y', X, X'$ is determined by the above definitions. An automorphism of $P(3)$ which restricts to the centre as $C \mapsto C^j$ sends $z$ to $jz$, and we may define $z$ to be the restriction from $\widetilde{P}$ of the generator of the same name. We also note that $\mathrm{Res}^{\widetilde{P}}_{P(3)}(c_2) = -xY' - x'Y - x^2 - x'^2$.

*Proof.* Much of the proof is exactly as in the case when $p > 3$, except that the expressions for $YY'$, $Y^2$ and $Y'^2$ require a different proof. Again we consider the spectral sequence for $\mathrm{B}P(3)$ as an $S^1$ bundle over $\widetilde{P}$, which has $E_2^{*,*} \cong \Lambda[u] \otimes H^*(\mathrm{B}\widetilde{P}; \mathbb{F}_p)$, and $d_2(u) = \pm yy'$. The proof that $x$, $x'$, $Y$ and $Y'$ form a basis for $H^2$ by calculating restrictions to the subgroups $\langle A, C \rangle$ and $\langle B, C \rangle$ works exactly as in Theorem 6. Many of the relations may be proved exactly as in Theorem 6, although the fact that $\langle y, y, y \rangle$ is equal to $x$ rather than 0 makes the relations look different. For example, $yY = y\langle y, y, y' \rangle \equiv \langle y, y, y \rangle y' = xy'$, where the congruence is modulo $\{0\}$. The relations containing $y'Y'$, $y'Y$, $yX$ and $y'X'$ follow similarly. As in the proof of Theorem 6, any $Z \in \langle (x, x'), \begin{pmatrix} y' \\ -y \end{pmatrix}, y \rangle$ satisfies $Zy = -xY + xx'$ and $Zy' = x'Y - xY'$. Similarly, any $Z' \in \langle (x', x), \begin{pmatrix} y \\ -y' \end{pmatrix}, y' \rangle$ satisfies $Z'y' = -x'Y' + xx'$ and $Z'y = xY' - x'Y$. We deduce that $X, X'$ satisfy the relations claimed in $H^4$, and that $X, X', xy, xy', x'y', Yy'$ form a basis for $H^3$. The relations $xY'y = x'^2 y$ and $x'Y'y = x^2 y'$ follow easily from the relations we have already proven, now we apply the Bockstein to them, and obtain $x(xY' + x'Y) = -xx'^2$ and $x'(xY' + x'Y) = -x'x^2$. It follows that the relation in $H^4$ yielding the spectral sequence relation $x[uy'] = x'[uy]$ must be $0 = xY' + x'Y + c_2 + x^2 + x'^2$. We deduce that $x^2, xx', x'^2, xY, xY', x'Y$, and $x'Y'$ form a basis for $H^4$. Now we shall return to the other relations we wish to prove in $H^4$.

The subspace of $H^4$ of elements restricting trivially to $\langle A, C \rangle$ and $\langle B, C \rangle$ contains $YY'$, and has basis $xx'$, $xY'$, $x'Y$. Hence there is some expression $YY' = \lambda xx' + \lambda' xY' + \lambda'' x'Y$. The automorphism $\theta$ of $P(3)$, which has the effect of exchanging the 'primed' and 'unprimed' generators fixes $YY'$, so $\lambda'$ and $\lambda''$ must be equal. Multiplying both sides of the expression for $YY'$ by $y$ gives $xx'y = \lambda xx'y - \lambda' x'^2 y$, and hence $\lambda = 1$, $\lambda' = 0$ as required. The expressions for $Y^2$ and $Y'^2$ may be deduced from that for $YY'$ by considering the effect of various automorphisms of $P(3)$.

We apply the Bockstein to the expressions for $Y^2$, $YY'$, $Y'^2$ and $\mathrm{Res}^{\widetilde{P}}_{P(3)}(c_2)$ respectively to obtain the relations $XY = -xX'$, $X'Y' = -x'X$, $XY' = -X'Y$ and $xX' = -x'X$. The spectral sequence argument suggested in Theorem 6 shows that $XY'$ is not in $xH^3 + x'H^3$, and so no new generators are needed in $H^5$. It is easily checked that $XX'$ is annihilated by both $x$ and $x'$, which implies that $XX' = 0$. We already have the relations $x^3 y' = x'^3 y$, $x^3 x' = x'^3 x$, and as in the case $p > 3$, we prove the relation $x^3 Y' + x'^3 Y = -x^2 x'^2$ by applying $\mathrm{P}^1$ to the relation $0 = x'Y + xY' + c_2 + x^2 + x'^2$, and we apply the Bockstein to this relation to obtain $x^3 X' + x'^3 X = 0$.

It may be checked that all products of degree at least 6 of elements of degree at most three may be expressed in the form

$$f_1 + f_2 Y + f_3 Y' \qquad \text{for even total degree}$$
$$f_1 y + f_2 y' + f_3 X + f_4 X' \qquad \text{for odd total degree}$$

where $f_i$ is a polynomial in $x$ and $x'$. The relations we have given between such elements are sufficient to imply the relations that hold between the corresponding elements of the spectral sequence, hence our presentation of the ring $H^*(\mathrm{B}P(3); \mathbb{F}_3)$ is complete. The effect of automorphisms on $z$ follows from its definition as a restriction from $\widetilde{P}$, and we already have the required expression for $\mathrm{Res}^{\widetilde{P}}_{P(3)}(c_2)$. ∎

It is unreasonable to expect that by choosing different generators we could make the statement of Theorem 7 look more like the statement of Theorem 6. For example, for $p > 3$, $H^2(\mathrm{B}P(3); \mathbb{F}_p)$ may be expressed as the direct sum of two $\mathrm{Aut}(P(3))$-invariant subspaces:
$$H^2(\mathrm{B}P(3); \mathbb{F}_p) = \langle x, x' \rangle \oplus \langle Y, Y' \rangle,$$
whereas for $p = 3$ the subspace $\langle x, x' \rangle$ has no $\mathrm{Aut}(P(3))$-invariant complement.

In each case, the action of the Steenrod algebra on the generators of $H^*(\mathrm{B}P(3); \mathbb{F}_p)$ is determined (up to a scalar multiple in the case of $\mathrm{P}^1(d_4)$) as in Proposition 3, except for the action of $\mathrm{P}^1$ on $X$ and $X'$. This is described below.

**Proposition 8.** *With notation as in Theorems 6 and 7,*
$$\mathrm{P}^1(X) = x^{p-1} X + zy$$
$$\mathrm{P}^1(X') = x'^{p-1} X' - zy'.$$

*Proof.* The spectral sequence operation ${}_B \mathrm{P}^1$ defined by Araki and Vasquez ([1], [16]) on the $E_\infty$ page of the spectral sequence for $\mathrm{B}P(3)$ as an $\mathrm{S}^1$-bundle over $\mathrm{B}\widetilde{P}$ sends $ux'$ to $ux'^p$, so we deduce that $P^1(X') \equiv x'^{p-1}X'$ modulo the image of the restriction from $\widetilde{P}$. Let $K$ be a subgroup of $P(3)$ of index $p$, and let $H^*(BK; \mathbb{F}_p) = \mathbb{F}_p[\bar{x}, c] \otimes \Lambda[\bar{y}, d]$, where $\beta(d) = c$, $\beta(\bar{y}) = \bar{x}$, and $d$, considered as a morphism from $K$ to $\mathbb{F}_p$, sends $C$ to 1. Then if we let $\text{Res} = \text{Res}^{P(3)}_K$ we have $\text{Res}(X') = \lambda(\bar{x}d - c\bar{y})$, since $X'$ is in the image of the Bockstein and the image of the Bockstein in $H^3(BK)$ is generated by $\bar{x}d - c\bar{y}$. We obtain

$$\text{Res}(P^1(X') - x'^{p-1}X') = -\lambda(c^p - \bar{x}^{p-1}c)\bar{y}.$$

If $P$ is an expression of degree $2p+1$ involving only $y, y', x$, and $x'$, then $\text{Res}(P)$ is a multiple of $\bar{x}^p \bar{y}$, and if $Q$ is in the span of $zy$ and $zy'$, then $\text{Res}(Q)$ is a multiple of $(c^p - \bar{x}^{p-1}c)\bar{y}$. We know that $P^1(X') - x'^{p-1}X' = P + Q$ for some such choices of $P$ and $Q$, and we deduce that for all $K$, $\text{Res}(P) = 0$. Thus $\text{Res}(\beta P) = 0$, and hence $\beta P$ is a multiple of $x^p x' - x'^p x$, so is zero in $H^*(BP(3); \mathbb{F}_p)$. The Bockstein is injective on the subspace of $H^{2p+1}$ generated by $x, x', y$, and $y'$, so we deduce that $P = 0$. Using explicit cochain calculations as in the proof of Theorem 6 we can determine $\lambda$ in the case $K = \langle A, C \rangle$ and $\bar{y} = \text{Res}(y)$ (resp. $K = \langle B, C \rangle$ and $\bar{y} = \text{Res}(y')$), and conclude that $Q = -zy'$. The result for $P^1(X)$ follows from this result or may be deduced similarly. ∎

**Acknowledgements.** The author thanks his research supervisor, C. B. Thomas, also D. J. Benson, D. J. Green and P. H. Kropholler. The author is indebted to J. F. Adams, who first suggested the use of Massey products in these calculations.